\documentclass[a4paper, 12pt]{article}
\usepackage{mathrsfs}
\usepackage{amsmath, amsthm, amsfonts, amssymb}
\usepackage{color}

\newtheorem{theorem}{Theorem}[section]
\newtheorem{theorem*}{Theorem*} 
\newtheorem{corollary}[theorem]{Corollary}
\newtheorem{lemma}[theorem]{Lemma}
\newtheorem{proposition}[theorem]{Proposition}
\theoremstyle{definition}
\newtheorem{definition}{Definition}[section]

\numberwithin{equation}{section}
 \theoremstyle{remark}
\newtheorem{remark}{Remark}[section]

\newcommand{\R}{\mathbb R}
\newcommand{\N}{\mathbb N}

\newcommand{\F}{\mathcal F}
\newcommand{\Pro}{\mathbb P}
\newcommand{\Epc}{\mathbb E}
\newcommand{\e}{\epsilon}

\begin{document}

\title{ 
Log-Harnack inequality  for reflected 
SPDEs driven by multiplicative noises and its applications}
\author{ Bin Xie\footnote{ E-mail:
 bxie@shinshu-u.ac.jp; bxieuniv@outlook.com}\\
{ \small{ Department of Mathematical Sciences, Faculty of
Science, Shinshu University }}\\
{ \small{3-1-1 Asahi, Matsumoto, Nagano
390-8621, Japan } } }
\date{  }
\maketitle


\begin{abstract}
We mainly investigate the log-Harnack inequality for 
the reflected stochastic partial differential equation driven 
by multiplicative noises  based on the  gradient estimate of the associated Markov semigroup.
To do it,  the penalization method and the comparison 
principle  for stochastic partial differential equations are adopted. As its applications, we  also establish the entropy-cost 
inequality and study some important estimates of the transition density relative to its invariant measure.
\end{abstract}



\section{Introduction}
Consider the following reflected stochastic partial differential equation (SPDE) on the bounded interval $[0,1]$ driven by a space-time white noise:
\begin{align}\label{spde-1.1}
\left\{
\begin{aligned}
\frac{\partial u}{\partial t}(t, x)&=\frac{1}{2} \frac{\partial^2 u}{ \partial x^2} (t,x) + b(u(t,x)) \\
&\qquad  + \sigma(u(t,x))\dot{W}(t,x) +\eta(dtdx), \ t>0, \ x\in (0, 1),\\
u(t,0)& =u(t,1)=0,\  t\geq 0,\\
u(0,x)& =h(x),\   x\in [0, 1],
\\
u(t,x) &  \geq 0,\  t>0, \ x\in [0, 1]\ a.s.,
\end{aligned}
\right.
\end{align}
where the  initial value $h$ is a non-negative, continuous function defined on the interval $[0,1]$ with $h(0)=h(1)=0$, the coefficients 
$b$ and $\sigma$ are measurable real-valued functions defined on $\R$,  $\{\dot{W}(t,x): t\geq 0, x\in [0,1] \}$ denotes the space-time white noise on a complete probability space 
$(\Omega, \mathcal{F},\Pro)$ with a filtration  $(\mathcal{F}_t)_{t\geq 0}$, and $\eta$ denotes a positive random measure on 
$[0, \infty) \times [0,1]$ and is supported on the zero set of $u$.
In this paper, without loss of generality, we consider the filtration  $(\mathcal{F}_t)_{t\geq 0}$ generated by $\dot{W}(t, x)$, that is, 
$\mathcal{F}_t:=\sigma(W(s,x): s\in [0, t], x\in [0,1])\vee \mathcal{N}$, where $\mathcal{N}$ denotes the family of $\Pro$-null sets.

The problem \eqref{spde-1.1} is one kind of random parabolic  obstacle 
problems and is also regarded as an infinite-dimensional Skorokhod 
problem. It is initially proposed  by D. Nualart and E. Pardoux 
\cite{NuPa-92} in 1992 for additive noises (i.e. $\sigma\equiv 1$) and 
then it is generalized to the case of multiplicative noises by C. 
Donati-Martin and E. Pardoux \cite{DoPa-93} in 1993.  About nine  
years later, T. Funaki and S. Olla \cite{FuOl-01} in 2001 proved that the 
fluctuation of a   $\nabla \phi$-interface model on a hard wall weakly 
converges  to the stationary solution of  the reflected SPDE driven by 
additive noises, which is known as a  famous and interesting 
application of reflected SPDEs. In addition,  various important 
properties, such as reversible measures 
\cite{FuOl-01}, \cite{Zam-01} and hitting 
properties \cite{DaMuZa-06},  of the solution to reflected SPDE 
\eqref{spde-1.1} driven by additive noises were studied.  

However, 
for the general case of the reflected SPDE \eqref{spde-1.1},
the study of it becomes more difficult. In fact, even for the 
uniqueness of its solution was left open in \cite{DoPa-93} for very 
long time and  at last  it was successfully solved in \cite{XuZh-09} until 2009. After that, the strong Feller property and the large deviation principle  of the  Freidlin-Wentzell type relative to  the solution were also investigated, see 
\cite{XuZh-09}, \cite{Zh-10} and references therein. 
For more information on the research of random obstacle problems, we also like to refer the readers to the monograph \cite{Zam-17} authored by L. Zambotti in 2017 and references therein.

The main purpose of this paper is to investigate the
 log-Harnarck inequality, a weak type of the dimension-free
  Harnack inequality, 
with respect to the Markov semigroup generated by the solution 
of \eqref{spde-1.1} by the approach of gradient estimate. The dimension-free Harnack inequality is
 initially introduced in \cite{Wan-97} by F.-Y. Wang to 
 study the log-Sobolev inequality of a diffusion process on
  Riemannian manifolds
and then it becomes as a very powerful and effective tool to 
the 
study  of various important properties of diffusion semigroups, such 
as, Li-Yau type heat kernel bound \cite{AeThWa-06}, 
hypercontractivity, ultracontractivity \cite{RoWa-03},  
strong Feller 
property, estimates on the heat kernels \cite{Wan-07} and Varadhan type small time asymptotics \cite{Zh-10}.
Generally speaking, the dimension-free Harnack inequality for the  Markov semigroup $(P_t)_{t\geq 0}$ relative to a Markov process with values in a Banach space $E$ is formulated as 
\begin{align}\label{eq-1.2-170521}
\psi(P_t\Phi(x)) \leq \big(P_t(\psi(\Phi))(y) \big) \exp\{\Psi(t,x,y) \}, \ x,y \in E, 0\leq  \Phi \in B_b(E), t>0,
\end{align}
where $\psi: [0, \infty) \to [0, \infty)$ is convex, $\Psi$ is a non-negative function defined on $[0, \infty) \times E \times E$ with 
$\Psi(t,x,x)=0$, and $B_b(E)$ denotes the family of all  measurable 
and bounded functions defined on $E$. In particular, for the special 
case $\psi(r)=r^p, p>1$, \eqref{eq-1.2-170521} 
is reduced to 
$$
(P_t\Phi(x))^p \leq \big(P_t(\Phi^p)(y) \big) \exp\{\Psi(t,x,y) \}, \ x,y \in E, 0\leq  \Phi \in B_b(E), t>0,
$$
which is called the Harnack inequality with power 
$p$ and for $\psi(r) =e^r,$ \eqref{eq-1.2-170521} is called the 
log-Harnack inequality, which, by (2) in Section 1.1.1 \cite{Wan-13},
 is equivalent to 
\begin{align*}
P_t\log\Phi(x) \leq \log P_t\Phi (y) + \Psi(t,x,y), \ x,y \in E, 0\leq  \Phi \in B_b(E), t>0.
\end{align*}
Here, we would also like to refer the readers to the monograph \cite{Wan-13} by F.-Y. Wang for more information of the  dimension-free Harnack inequality. 

Recently, such dimension-free Harnack inequalities have been 
actively studied and applied to the field of stochastic partial differential equations.  
For example, we refer the readers  to \cite{DaRoWa-09}, 
\cite{Liu-09}, 
\cite{Wan-07}, \cite{Wan-15}, \cite{Zh-10} for  the study of the Harnack inequality with power, and respectively to 
\cite{RoWa-10}, \cite{WaWuXu-11}, \cite{WaXu-14}, \cite{WaZh-14}  for the study of the log-Harnack inequality.  
In particular, in \cite{Zh-10}, T.-S. Zheng studied the  Harnack 
inequality with power relative to the solution of \eqref{spde-1.1} 
driven by additive noises by the coupling method with  the help of the Girsanov theorem and then applied it to the
study of strong Feller property, hyperbound property and others. 

However, most of the works listed above studied the SPDEs with additive (colored or white) noise and monotonic drifts.
The studies on the Harnack inequalities with power for  semilinear stochastic partial differential equations driven by  multiplicative noises are very rare. It is initially established by F.-Y. Wang in \cite{Wan-11} and then is generalized in \cite{ZhSa-13}, see also \cite{Wan-15}. 
To the best of my knowledge, all of the works deal with a noise ``multiplicative noises'', roughly speaking,  the diffusion coefficients are Hilbert-Schmidt operator even there is no reflection, see 
\cite{RoWa-10} or \cite{Wan-15} for more 
explanations. On the other hand, in our case, the space-time white noise is considered, which is more difficult.
So, under our framework, it seems to be impossible to establish  the  Harnack inequality at present, because the space-time white noise is considered and the diffusion coefficient $\sigma$ can not be  regarded as a Hilbert-Schmidt operator on $L^2(0,1)$.  Therefore, in this paper,  instead of the Harnack inequality with power, in our framework, we devote to investigating the weak form of the Harnack inequality, i.e.,  the log-Harnack inequality to  the solution of \eqref{spde-1.1}. 
The log-Harnack inequality relative to SPDEs driven by multiplicative 
noises is first studied in \cite{RoWa-10} for a special  case of the 
coefficient $\sigma$  and then it is generalized in  \cite{WaZh-14} 
to the general coefficient $\sigma$. All of these 
works deal with  the  SPDEs  without reflection term $\eta$. It is 
valuable to point out that  their results obtained in  \cite{RoWa-10}  
and \cite{WaZh-14} can not be applied to reflected SPDEs 
\eqref{spde-1.1}, see Remark \ref{rem-4.1-170607} for details.
This is one of the main motivations of this paper. To investigate the 
log-Harnack inequality relative to the solution $u(t,x)$ of 
\eqref{spde-1.1},  
the estimate of the Gradient estimate of the solution is studied by 
the approach of  penalization and a comparison principle
 of SPDEs.
 
Let us introduced some notations, which will be used throughout
 this paper. Let $H=L^2(0,1)$ and denote its scalar inner 
 product and norm 
by $\langle \cdot, \cdot \rangle$ and $|\cdot|$ respectively, that is 
$$
\langle h, \tilde{h} \rangle=\int_0^1 h(x) \tilde{h}(x)dx\  \ 
\mbox{and} \  \ |h|=\langle h, h \rangle^{\frac12}, \ h, \tilde{h} \in H.
$$
Let $K_0:=\{h\in H: h(x) \geq 0, x\in [0,1]\}$, the family of all non-
negative functions $h\in H$.   It is known that $K_0$  is the closed 
subset set of $H$, and in particular it is Polish. Furthermore, 
we denote by $C_0^k(0,1)$  the class of all $k$-th continuously differentiable functions $h$ defined $[0,1]$ with $h(0)=h(1)=0$. For
simplicity, we write $C_0(0,1)$ for $C_0^0(0,1).$

Let us now formulate the definition of the reflected SPDE \eqref{spde-1.1} according to the initial introduction in \cite{DoPa-93}.
\begin{definition} Suppose the initial value $h \in K_0\cap C_0(0,1)$.
A pair $(u,\eta)$ is said to be a solution of the reflected SPDE \eqref{spde-1.1} if the following are satisfied.\\
{\rm(i)} The process $\{u(t,x): t \geq 0, x\in [0,1]\}$ is non-negative and continuous on 
$[0, \infty)\times [0,1]$. Moreover, $(u(t, \cdot))_{t\geq 0}$ is $\F_t$-adapted.  \\
{\rm(ii)} $\eta$ is a positive random measure on $[0, \infty) \times [0,1]$ such that\\
\quad {\rm(a)} $\eta(\{t\} \times (0,1))=0,\ t\geq 0;$\\
\quad {\rm(b)} $\int_0^t\int_0^1 x(1-x) \eta(dsdx)<\infty,  \ t\geq 0;$\\
\quad {\rm(c)} $\eta$ is $\F_t$-adapted, that is, for any Borel set $B\subset [0,t] \times [0,1]$, $\eta(B)$ is $\F_t$-\\
\qquad measurable.\\
{\rm(iii)} For any $\phi\in C_0^2(0,1)$,  the following stochastic integral equation is 
fulfilled.
\begin{align*}
\langle u(t), \phi\rangle =& \langle h, \phi\rangle +\frac{1}{2}\int_0^t \langle u(s), \phi'' \rangle ds + \int_0^t \langle b(u(s)), \phi \rangle ds\\
& 
+ \int_0^t\int_0^1 \sigma(u(s,x))\phi(x)W(dsdx) +\int_0^t\int_0^1 \phi(x)\eta(dsdx), \notag\  t\geq 0 \ a.s.,
\end{align*} 
 where $u(t):=u(t,\cdot)$ and the first term in second line is understood in the sense 
of It\^{o}'s integral, see \cite{Wal}.\\
{\rm (iv)} The support of $\eta$ is a subset of $\{(t,x): u(t,x)=0, t\geq 0, x\in [0,1]\}$, that is, $$\int_0^\infty \int_0^1 u(t,x)\eta(dtdx)=0\ a.s.$$
\end{definition}
\begin{remark}
It is clear that {\rm(c)} in  {\rm(ii)} is equivalent to the following:\\
the stochastic process
$ \int_0^t\int_0^1 \phi(s,x)\eta(dsdx), t \geq 0$
is $\F_t$-measurable for any measurable function $\phi: [0, \infty)\times [0,1] \to [0, \infty)$. 
\end{remark}

Based on \cite{DoPa-93}, let us formulate the existence and uniqueness of the solution of \eqref{spde-1.1} under the following assumptions on  $b$ and $\sigma$.
\\
{\bf A1)} $b: \R \to \R$ is globally Lipschitz continuous, that is, there exists a 
constant $L_b>0$ such that for any $u, v\in \R$
$$
|b(u) - b(v)| \leq L_b |u-v|.
$$
{\bf A2)} $\sigma: \R \to \R$ is globally Lipschitz continuous, that is, there exists a 
constant $L_\sigma>0$ such that for any $u, v\in \R$
$$
|\sigma(u) - \sigma(v)| \leq L_\sigma |u-v|.
$$
To clarify the relations between our results and the Lipschitz constants of $b$ and $\sigma$ as below, we described the above assumptions separately. 

Let us now formulate the result of the existence and uniqueness of the solution  to the reflected SPDE \eqref{spde-1.1}. 
\begin{theorem}\label{thm-1.1-170501}
Suppose  the assumptions {\bf A1)} and {\bf A2)} are satisfied. Then for each $h\in K_0\cap C_0(0,1)$, the reflected SPDE \eqref{spde-1.1} has a unique solution $(u, \eta)$  such that for any $T>0$ and $p \geq 1$
$$ \Epc \left[\sup_{(t,x)\in [0,T] \times [0,1]}|u(t,x)|^p \right]<\infty.$$
\end{theorem}

\begin{remark}
{\rm (i)} 
We refer the readers to Theorem 3.1 \cite{DoPa-93} and Theorem 2.1 \cite{XuZh-09} for the proof of Theorem \ref{thm-1.1-170501}.
For the case of  multiplicative noises, the  uniqueness of the solution 
is more difficult. In fact, although in 1993, C. Donati-Martin and E. 
Pardoux \cite{DoPa-93} succeeded in the  construction a solution of   
\eqref{spde-1.1} by the approach of penalization and showed that it is 
minimal, the uniqueness of the solution  is left as an open problem.  
It is  at last solved by T.-G. Xu  and T.-S. Zhang 
\cite{XuZh-09} until 2009 . \\
{\rm (ii)} The H\"{o}lder continuity of $u(t,x)$ is also studied in 
\cite{DaZh-13}.  It is  proved that $u(t,x)$ is H\"{o}lder continuous 
with  the index $(\frac14-, \frac12-)$ on $[0,T] \times [0,1]$, which 
is same as that of the solution to its penalized equation, see 
\eqref{spde-approx} below.
\end{remark}

Throughout this paper, we will sometimes make use of  the notations 
$u(t)$, $u(t,x;h)$ or $u(t;h)$ instead of $u(t,x)$ according to our purposes. 
For example, whenever we are going to emphasize its initial value $h$ and regard the solution $u(t,x)$ as an $H$-valued process, we will use the notation $u(t;h)$.

This paper is organized as follows: 
In Section \ref{sec-2}, the main results, including the 
gradient estimate and the log-Harnack inequality for the solution of \eqref{spde-1.1},
are formulated. 
In addition, some applications of  the log-Harnack inequality are stated with brief proofs.
From Section  \ref{sec-3}, we devote to the proofs of 
the main results.  In Section \ref{sec-3}, we give a brief proof of 
Proposition  \ref{thm-2.1-170524}, which shows that the solution is 
continuous in its initial value. The proofs of the gradient 
estimate and the log-Harnack inequality relative to the solution are stated 
in Section \ref{sec-4} by  the penalization method and  a comparison theorem principle, and then as another  application of the gradient estimate,  the proof of Theorem  \ref{thm-2.4-170605} is described in last section.

\section{Main Results}\label{sec-2}
In this section, we first state the continuity of the solution $u(t;h)$ 
on its initial data. 
After that, we formulate the main results of this paper including  the gradient estimate and the log-Harnack inequality for the solution $u(t;h)$ of \eqref{spde-1.1} and finally state some  applications of the log-Harnack inequality with  brief proofs.
 
\begin{proposition}\label{thm-2.1-170524}
Suppose assumptions in Theorem \ref{thm-1.1-170501} are satisfied and let the pair $(u(t; h_1), \eta_1)$ and respectively $(u(t;h_2), \eta_2)$ denote the unique solutions of \eqref{spde-1.1} with initial data $h_1$ and $h_2 \in K_0\cap C_0(0,1)$. Then
for each $p\geq 1$ and $T>0$,  there exists a constant $K=K(T,p)>0$ such that 
\begin{align}\label{eq-2.1-170524}
\Epc \left[\sup_{(t,x)\in [0,T] \times [0,1]}|u(t, x; h_1)- u(t, x; h_2)|^p \right]\leq K|h_1-h_2|_{\infty}^{p},
\end{align}
where $|h|_{\infty}= \sup_{x\in [0,1]} |h(x)|$ for $h\in C_0(0,1)$.
\end{proposition}

From now on, we will formulate the main theorems of this paper.
Let us here introduce some notations. Let $B_b(H)$ (respectively 
$C_b(H)$) denote the class of all bounded and measurable real-valued 
functions (respectively bounded and continuous real-valued 
functions) defined on $H$. Similarly, we will use the notations $B_b(K_0)$ and $C_b(K_0)$. 
Let us denote $\Pi: H \to K_0$ the projection  from 
$H$ to the closed convex set $K_0$. Then, we can identify 
 $B_b(K_0)$ (respectively $C_b(K_0)$) with a subspace of 
 $B_b(H)$ (respectively $C_b(H)$) by means of the injection:
  $B_b(K_0)\ni \Phi \mapsto \Phi(\Pi) \in B_b(H)$, see Section 2 \cite{Zam-01} for details.

Let us define $P_t$ by
$$
P_t\Phi(h)= \Epc [\Phi(u(t;h))], \  t\geq 0, h\in B_b(K_0).
$$
It is shown that $(P_t)_{t\geq 0}$ forms a Markov semigroup \cite{Zh-10}.
To study the gradient estimate and the log-Harnack inequality  relative to $P_t$, we will further suppose the following assumption {\bf A3)}.\\
{\bf A3)} There exist two constants $0< \kappa_1< \kappa_2$, such that for $u\in \R$, 
$$ 
\kappa_1 \leq |\sigma(u)| \leq \kappa_2.
$$
Throughout this paper, for any function $\Phi$ defined on $H$, we denote by $|\nabla \Phi|(h)$ the local Lipschitz constant of $\Phi$
at $h\in H$, i.e., 
$$
|\nabla \Phi|(h) =\limsup_{|\tilde{h} -h|\to 0} \frac{|\Phi(\tilde{h}) -\Phi(h)|}{|\tilde{h} -h|}. 
$$ 
To state the main results, let us introduce 
the constant $M=M(L_b, L_\sigma)>0$ and the function $\zeta(t), t\geq 0$ as below. 
\begin{align}\label{eq-3.16-170515}
M&= \max \left\{ 
3, \frac{9L_\sigma^2}{\sqrt{\pi} }, \frac{8L_b^2}{L_\sigma^{4}}, \frac{144 L_b^2}{ L_\sigma^{2}\sqrt{\pi}}, \frac{864L_b^2}{\sqrt{\pi}} \right\},
\end{align}
and
\begin{align}
\zeta(t) & =t^{\frac{1}{2}}+ \frac{9L_\sigma^4}{4}t + \frac{3L_b^2}{2}t^2 + \frac{18L_b^2L_\sigma^2}{5\sqrt{\pi}} t^{\frac52},\  t\geq 0. \label{eq-2.3-170604} 
\end{align}
It is important to see that the constant  $M$ and the coefficients of $\zeta(t), t\geq 0$ depend only on $L_b$ and $L_\sigma$. The first main theorem is the gradient estimate of $P_t$, which plays a key role to the log-Harnack inequality.
\begin{theorem}\label{thm-2.2-170519}
Suppose the assumptions ${\bf A1)}$-{\bf A3)} are satisfied. Then we have that  for any $\Phi\in C_b^1(K_0)$, 
\begin{align}\label{eq-2.4-170611}
|\nabla P_t \Phi|^2 \leq 2M\exp\{\zeta (t) \} P_t |\nabla  \Phi|^2,
 \ t \geq 0.
\end{align}
\end{theorem}
According to the approaches used in \cite{RoWa-10} and \cite{WaZh-14}, the log-Harnack inequality of $P_t$ can be obtained based on Theorem  \ref{thm-2.2-170519}.
\begin{theorem}\label{thm-2.3-170519}
Under the assumptions of Theorem  \ref{thm-2.2-170519},  the log-Harnack inequality  holds for $P_t, t>0$. More precisely,   for any  strictly positive $\Phi\in B_b(K_0)$, 
\begin{align}\label{eq-2.4-170520}
P_t \log \Phi(h_1) \leq \log P_t\Phi(h_2) +  \frac{M|h_1-h_2|^2}{\kappa_1^2\int_0^t \exp\{-\zeta(s) \}ds}, \ h_1, h_2\in K_0, t>0, 
\end{align}
where $\kappa_1$ is the constant in {\bf A3)}.
\end{theorem}
\begin{remark}
Because of  the multiplicative noise is considered, we just establishes the log-Harnack inequality as above. However, for the case of the additive noise, the Harnack inequality can be established, see \cite{Zh-10} for the Lipschitz case and \cite{Xie-18} for the nonlinear case.
\end{remark}

By Theorem  \ref{thm-2.2-170519}, we can  also obtain    a estimate for $P_t\Phi^2 - (P_t\Phi)^2$ and the Lipschitz continuity of $P_t \Phi$ as below. 
\begin{theorem}
\label{thm-2.4-170605}
Under the assumptions of Theorem \ref{thm-2.2-170519}, we have that\\
{\rm (i)} 
For any $\Phi\in C_b^1(K_0)$, 
\begin{align}\label{eq-2.7-170611}
P_t\Phi^2 - (P_t\Phi)^2 \leq 2Mk_2^2 P_t |\nabla  \Phi|^2  \int_0^t \exp\{\zeta(s) \}ds, \ t \geq 0,
\end{align}
where $\kappa_2$ is the constant in  {\bf A3)}.
\\
{\rm (ii)}  
For any  $\Phi\in B_b(K_0)$, 
\begin{align}\label{eq-2.6-170611}
|\nabla P_t \Phi|^2 \leq  \frac{2M}{k_1^2\int_0^t \exp\{-\zeta(s) \}ds}  \{P_t\Phi^2 - (P_t\Phi)^2\}, \ t > 0.
\end{align}
In particular, $P_t \Phi$ is Lipschitz continuous for  each $\Phi\in B_b(K_0)$ and $t>0$.
\end{theorem}

From now on, let us formulate some applications of the log-Harnack inequality of $P_t$. Let us first study the  entropy-cost inequality. 
To formulate it, let us introduce the definition of $L^2$-Wasserstein distance between two probability measures.
For any two probability measures $\nu$ and $\mu$ on $K_0$, we denote by
$W_2(\nu, \mu)$ the $L^2$-Wasserstein distance between them with respect to the cost function $(h_1, h_2) \to |h_1-h_2|$, i.e.,
$$
W_2(\nu, \mu)=\inf_{ \mathcal{P} \in \mathcal{C}(\nu, \mu)}
\left\{ 
\int_{K_0\times K_0} |h_1-h_2|^2 \mathcal{P} (dh_1, dh_2)
\right\}^{\frac12},
$$
where $\mathcal{C}(\nu, \mu)$ denotes the family of all couplings
of  $\nu$ and $\mu$. 
It is also called $L^2$-transportation cost between $\nu$ and $\mu$.

Under the assumptions {\bf A1)}-{\bf A3)},  the existence of  the
invariant probability measure of the  reflected SPDE \eqref{spde-1.1} can be shown by imitating 
the idea introduced in Theorem 
2.1 \cite{ZhYa-14}, in which stochastic heat equation with double reflections is studied.   See also  \cite{Kal-18} for a proof in our case. In the following,  we denote the invariant measure by  $\pi$ and regard it as the invariant measure of $P_t$ on $K_0$.
As the application of Theorem \ref{thm-2.3-170519}, we first have the following entropy-cost inequality for the adjoint operator of $P_t$ in 
$L^2(\pi)$.
\begin{corollary}
Let $\pi$ be one invariant measure of $P_t$  and let $P_t^*$ be the adjoint operator of $P_t$ in $L^2(\pi)$. Then for any non-negative $\Phi\in B_b(K_0)$  with $\pi(\Phi)=1$, we have the following entropy-cost inequality:
\begin{align*}
\pi\big((P_t^*\Phi) \log P_t^*\Phi \big) \leq \frac{MW_2(\Phi\pi, \pi)^2}{\kappa_1^2\int_0^t \exp\{-\zeta(s) \}ds}, \ t>0,
\end{align*}
where $\kappa_1$ is the constant in {\bf A3)}.
\end{corollary} 
The relation between  the log-Harnack inequality and the  entropy-cost inequality is initially observed in Corollary 1.2 
\cite{RoWa-10} and then it is described for different frameworks 
\cite{WaWuXu-11} and \cite{WaZh-14}. This corollary can be easily shown by analogy from Corollary 1.2 \cite{RoWa-10}.
Simply speaking, it is enough to apply \eqref{eq-2.4-170520} to 
$P_t^*\Phi$ in stead of $\Phi$ and then make use of the invariance of $\pi$. The
details of the proof are omitted here.

Finally, let us apply the log-Harnack inequality to the
entropy inequalities relative to the transition density 
of  the transition semigroup  $P_t$ with respect to its invariant measure $\pi$. 
\begin{corollary}
The Markov semigroup $(P_t)_{t\geq 0}$ is strong Feller, that is, for any $\Phi \in B_b(K_0)$, 
\begin{align*}
\lim_{|\tilde{h}-h|\to 0} P_t\Phi(\tilde{h}) =P_t\Phi(h),\ h\in K_0, t>0,
\end{align*}
and  
$P_t$ is irreducible, i.e., for any open set 
$\emptyset \neq A \subset K_0$, $P_t1_A(h)>0 $ for all $h\in K_0$. 
In particular, 
the invariant measure of $P_t$ is unique and is fully supported on $K_0$, which will be still denoted by $\pi$ in the following.
 
Moreover, for each $t>0$, $P_t$ is absolutely continuous 
with respect to its unique invariant measure $\pi$ with a strictly positive transition density $p_t( \tilde{h}, h)$ and the following estimates relative to the transition density $p_t(\tilde{h}, h)$ are satisfied.\\
{\rm (i)} The log-Harnack inequality \eqref{eq-2.4-170520} is equivalent to  the heat kernel  entropy inequality
\begin{align*}
\int_{K_0} p_t(h_1, h)\log\frac{p_t(h_1, h)}{p_t(h_2, h)}\pi(dh) \leq  \frac{M|h_1-h_2|^2}{k_1^2\int_0^t \exp\{-\zeta(s) \}ds},\ h_1, h_2\in K_0,   t>0. 
\end{align*}
{\rm (ii)} For any $h_1\in K_0$, 
\begin{align*}
& \int_{K_0} p_t(h_1, h) \log p_t(h_1, h)  \pi(dh) \\
\leq & -\log \int_{K_0} \exp \left\{  
\frac{-M|h_1-h|^2}{k_1^2 \int_0^t \exp\{-\zeta(s) \}ds }  \right\} \pi(dh),\ t>0.
\end{align*}
{\rm (iii)} For any $h_1, h_2\in K_0$,
\begin{align*}
 \int_{K_0} p_t(h_1, h) p_t(h_2, h)d\pi(h) \geq \exp\left\{ -\frac{M|h_1-h_2|^2}{k_1^2\int_0^t \exp\{-\zeta(s) \}ds} \right\}, \ t>0. 
\end{align*}
\end{corollary}
\begin{proof}
The above results are well-known as the applications of the log-Harnack inequality.  
Let us only state the proof of the irreducibility of $P_t$ and  refer the readers to Theorem 1.4.2 (2) \cite{Wan-13} for the proof of {\rm (i)} and for the others,  to Theorem 1.4.1 \cite{Wan-13} for the abstract results and to \cite{RoWa-10}, \cite{WaWuXu-11},  
\cite{WaZh-14} for different frameworks. For irreducibility, it is enough to show that 
$P_t1_A(h)>0, \ t>0$
with 
$A=A(\tilde{h}):=\{h: |\tilde{h} -h| < \gamma \}$ for any $\tilde{h}$ and 
$\gamma>0$. 
By the continuity of $u(t,x)$ on $[0, \infty)\times [0,1]$, for each $\gamma >0$, there exists $\tilde{t}>0$ such that    
$P_t1_A(\tilde{h}) =\Pro(|u(t; h) -h|< \gamma) >\frac12$.
By contradiction, we assume there exists $h_0\in K_0$ such that 
$P_t1_A(h_0)=0$. Applying \eqref{eq-2.4-170520} to $1 +n1_A, n \in \N$, we have
$\log P_t(1+ n1_A)(h_0)=0$ and then 
 for $t<\tilde{t}$ and all $n\in \N$,
\begin{align*}
\frac12\log(1+n) \leq P_t \log (1+n1_A(\tilde{h})) \leq  \frac{M|\tilde{h}-h_0|^2}{\kappa_1^2\int_0^t \exp\{-\zeta(s) \}ds},
\end{align*}
which is impossible,  because of the boundedness of the last term. Therefore,  for any $t\leq \tilde{t}$, 
$P_t$ is irreducible  and then for all $t>0$ by the Markov property of $P_t$.
\end{proof}


\section{Proof of Proposition \ref{thm-2.1-170524}}\label{sec-3}
Let us first summarize the following lemma according to Theorem 4.1 \cite{NuPa-92} as preparation. 
Let $\mathcal{C}_+$ be the family of  continuous functions $v(t,x)$  defined on $[0, \infty)\times [0,1]$ such that 
$v(t,0)=v(t,1)=0, t\geq 0$ and  $v(0,x)=h(x)$ such that 
$h(x)\in K_0 \cap C_0(0,1)$. 
\begin{lemma} \label{lem-3.1-170524}
$(1)$ For each $v\in \mathcal{C}_+$, there exists a unique pair $(z,\eta)$
such that the following are satisfied.
\\
{\rm(i)} $z(t,x)$ is continuous on $[0, \infty) \times [0,1]$  such that $z(0,x)=0, x\in [0,1]$, $z(t, 0)=$
\\
\qquad  $z(t,1)=0$ and $z(t,x)+v(t,x) \geq 0.$\\
{\rm(ii)} $\eta$ is a positive measure on $[0, \infty)\times [0,1]$ with $\eta([0,T] \times(\e, 1-\e ))<\infty$ 
for all 
\\
\qquad   $\e \in (0, \frac12)$ and each  $T>0$.\\
{\rm (iii)} For any $t\geq 0$ and $\phi\in C_0^2(0,1)$,
\begin{align*}
\langle z(t), \phi \rangle = \frac12 \int_0^t \langle z(s), \phi'' \rangle ds +\int_0^t\int_0^1 \phi(x)\eta(dsdx).
\end{align*}
{\rm (iv)} $\int_0^\infty \int_0^1 \big(z(t,x)+v(t,x) \big)\eta(dtdx)=0.$\\
$(2)$ Let $v_i \in \mathcal{C}_+$ and $(z_i, \eta_i)$ be the unique pair corresponding  to  
$v_i, \ i=1,2$ as above. Then for each fixed $T>0$ 
\begin{align*}
|z_1-z_2|_{T, \infty} \leq |v_1-v_2|_{T, \infty}.
\end{align*}
Hereafter, $|z|_{t, \infty} :=\sup_{(s,x)\in [0,t] \times [0,1]}|z(s,x)|, t>0$ for a continuous function $z(t,x)$ defined on $[0, \infty) \times [0,1]$.
\end{lemma}
\begin{remark}
The first part $(1)$ in  Lemma \ref{lem-3.1-170524} is the content of Theorem 4.1 \cite{NuPa-92} and the second part $(2)$ can been 
deduced from $(13)$ in Section 2 \cite{NuPa-92}.
\end{remark}

Using the above lemma, let us describe the proof of Proposition  \ref{thm-2.1-170524} briefly.\\
\noindent {\bf Proof of Proposition  \ref{thm-2.1-170524}:} 
To show this proposition, it is enough to prove that there exists a $p_0>1$ such that for any $p>p_0$ 
\begin{align}\label{eq-3.2-170525}
\Epc \left[ \left|u(t; h_1)- u(t; h_2) \right|_{T, \infty}^p \right]\leq
 C|h_1-h_2|_\infty^p
\end{align}
holds for a constant $C=C(T, p)>0$.
In fact, if this is true, then \eqref{eq-2.1-170524} can  be easily proved by Jensen's inequality.

Let $v_i(t,x), i=1,2$ be defined by 
\begin{align*}
v_i(t,x)= &\int_0^1g(t,x, y)h_i(y)dy +\int_0^t \int_0^1 g(t-s,x,y)b(u_i(s,y;h_i))dsdy\\
&+\int_0^t\int_0^1g(t-s, x, y)\sigma(u_i(s, y;h_i))W(dsdy). \notag
\end{align*} 
Hereafter, the function $g(t,x, y)$ denotes  the fundamental solution of the linear part of \eqref{spde-1.1}.
Then it is well-known that $v_i(t,x), i=1,2$ satisfy all of the assumptions in Lemma \ref{lem-3.1-170524}. 
Set
$$
z_i(t,x)= u_i(t,x; h_i)- v_i(t,x), \ t\geq 0, x\in [0,1], i =1,2.
$$
From the proof of Theorem 2.1 \cite{XuZh-09},
it follows that $(z_i, \eta_i)$ is the unique pair which 
satisfies all the conditions {\rm (i)}-{\rm (iv)} in Lemma \ref{lem-3.1-170524} for $i=1,2$.
Hence by Lemma \ref{lem-3.1-170524}, we have 
\begin{align*}
|z_1-z_2|_{T, \infty} \leq  |v_1-v_2|_{T, \infty}
\end{align*} 
and then
\begin{align*}
|u_1-u_2|_{T, \infty} \leq  2|v_1-v_2|_{T, \infty},
\end{align*}
which in particular implies that
\begin{align}\label{eq-3.3-170525}
\Epc [|u_1-u_2|_{T, \infty}^p]\leq 2^p \Epc[|v_1-v_2|_{T, \infty}^p], \ p>1.
\end{align} 
From the definitions of $v_i, i=1,2$,  it follows that for any $p>1$
\begin{align*} 
& \Epc \left[|v_1-v_2|_{T, \infty}^p \right] \\
\leq & 3^{p-1}\sup_{(t,x)\in [0, \infty)\times [0,1]}\left|\int_0^1 
g(t,x,y)(h_1(y)-h_2(y)) dy\right|^p \\
&+  3^{p-1} \Epc \left[\left|\int_0^\cdot \int_0^1 g(\cdot -s, \cdot,y) \Big(b(u_1(s, y; h_1))-b(u_2(s,y;h_2)) \Big)dsdy
\right|_{T, \infty}^p\right] \notag  \\
& 
+ 3^{p-1} \Epc \left[\left|\int_0^\cdot \int_0^1 g(\cdot -s, \cdot,y) \Big(\sigma(u_1(s, y; h_1))-\sigma(u_2(s,y;h_2)) \Big)W(dsdy)\right|_{T, \infty}^p \right] \notag \\
=: & 3^{p-1} \big(I_1+ I_2(T) +I_3(T) \big). \notag
\end{align*}
By the property of $g(t,x,y)$, it is easy to know that 
\begin{align*}
I_1\leq |h_1 - h_2 |_\infty^p, \ t\geq 0.
\end{align*}
Under the assumptions {\bf A1)}-{\bf A2)}, using the H\"{o}lder inequality and the Burkholder inequality, we have that for each $p>20$, there exists a constant $C_1=C_1(T, p)>0$ such that 
\begin{align*}
I_2(T) +I_3(T) \leq C_1(L_b^p+L_\sigma^p) \int_0^T 
\Epc [ |u_1(\cdot ;h_1) - u_2(\cdot;h_2)|_{t, \infty}^p]dt;
\end{align*}
see the proof of Theorem 2.1 \cite{XuZh-09}.
Combining the above estimates  with \eqref{eq-3.3-170525}, we have that there exists a constant $C_2=C_2(T, p, L_b, L_\sigma)>0$ such that 
 \begin{align*}
\Epc [|u_1-u_2|_{T, \infty}^p] \leq   C_2 |h_1 - h_2 |_\infty^p 
 +C_2 \int_0^T \Epc[ |u_1(\cdot ;h_1) - u_2(\cdot;h_2)|_{t, \infty}^p]dt. 
\end{align*}
Consequently, using Gronwall's inequality, we can conclude the proof of \eqref{eq-3.2-170525}, which completes  the proof. 
\hfill $\Box$

\section{Proofs of Theorem \ref{thm-2.2-170519} and Theorem \ref{thm-2.3-170519} } \label{sec-4}

We will adopt the approach of penalization initially introduced in 
\cite{DoPa-93} for our purpose.  The main idea to prove Theorem 
\ref{thm-2.2-170519} is based on the comparison principle for SPDEs. 

Let $e_n(x)=\sqrt{2}\sin(n\pi x), x\in [0,1]$.  Then it is easy to see that 
$$\frac{d^2}{dx^2}e_n(x)= -n^2\pi^2e_n(x),\ n\in \N$$ and  $\{e_n(x)\}_{n\in \N}$ forms 
a complete orthonormal basis of $H$.
\\ Set
$$
D(A)=\left\{h\in H: \sum_{n=1}^\infty n^4 \langle h, e_n\rangle^2 < \infty \right\}
$$
and
$$Ah=- \frac{1}{2}\sum_{n=1}^\infty n^2\pi^2  \langle h, e_n\rangle e_n, \ h\in D(A).$$
Then the operator $A$ with its domain $D(A)$  is the closure in $H$ 
of $\frac{1}{2}\frac{d^2}{dx^2}$ on $C_0^2(0,1)$. Remark that $D(A)$ 
coincides with the Sobolev space $H^2(0,1) \cap H_0^1(0,1)$.
\\
Let
$$
B(u)(x)=b(u(x)), \ x\in [0,1], u\in H,$$ 
and 
$$ \Sigma(u)[h](x)=\sigma(u(x))h(x), \  x\in [0,1], u, h\in H.
$$
Then, {\bf A1)}  gives that $B$ is a Lipschitz continuous mapping from 
$H$ to $H$. By  {\bf A3)}, we know that  $ \Sigma:H \to \mathcal{L}(H)$, 
where $\mathcal{L}(H)$ denotes the class of all bounded linear operators on $H$. 
However,  we point out that {\bf A2)} does not implies the  Lipschitz continuity of $\Sigma$ on $H$, which is different from $B$.
Define $w_n(t)$ by
$$
w_n(t) =\int_0^t\int_0^1 e_n(x)W(dsdx),\ t\geq 0.
$$
It is easy to know that $\{w_n(t), t \geq 0\}_{n\in \N}$ is a sequence of independent standard one-dimensional Brownian motions.
Set $$
W(t) =\sum_{n=1}^\infty w_n(t)e_n, \ t\geq 0.
$$
Then, the stochastic process $(W(t))_{t\geq 0}$ is a cylindrical 
Wiener process on $H$, and moreover we have that 
$$
\int_0^t\int_0^1 \sigma(u(s,x))W(dsdx)=\int_0^t \Sigma(u(s))dW(s),
$$
where the stochastic integral in the right hand side is understood in the sense of It\^{o}'s
integral with respect to a cylindrical Brownian motion on $H$, see \cite{DaZa-92} or \cite{DaZa-96}.

Let  $f:\R\to [0, \infty)$ be a non-increasing and Lipschitz continuous   function such that 
$$
f(u) \equiv 0, u\geq 0 \ \mbox{and}\ f(u)>0, u<0.
$$
For example, we can take $f(u)=u^{-}:= \max\{-u, 0\}$ or $f(u)= \arctan((u\wedge 0)^2)$. In particular, $f(u)= \arctan((u\wedge 0)^2)$ is twice continuously differentiable and bounded.

Let us consider the following penalized SPDE, which is used to construct the minimal solution of \eqref{spde-1.1} in \cite{DoPa-93}. For each $\e>0$ and $h\in H$, let us  consider the SPDE
\begin{align}\label{spde-approx}
\left\{
\begin{aligned}
\frac{\partial u^\e}{\partial t}(t, x)&=\frac{1}{2} \frac{\partial^2 u^\e}{ \partial x^2} (t,x) + b(u^\e(t,x)) +\frac{1}{\e}f(u^\e(t,x)) \\
&\qquad  + \sigma(u^\e(t,x))\dot{W}(t,x), \ t>0, \ x\in (0, 1),\\
u^\e(t,0)& =u^\e(t,1)=0,\  t\geq 0,\\
u^\e(0,x)& =h(x),\   x\in(0, 1).\\
\end{aligned}
\right.
\end{align}
In the following to emphasize the coefficients $b, f$ and $\sigma$, we will sometimes denote this  equation by SPDE$(b,f, \sigma)$.
Under our assumptions, it is known that for each $\e>0$, \eqref{spde-approx} has a unique mild solution $u^\e(t,x; h)$, which forms a Markov process on $H$. 

Moreover, it is proved that 
$u^\e(t,x) \geq u^{\e'}(t,x)\ a.s.$ for $\e< \e'$ and for any $p\geq 1$,
\begin{align*}
\lim_{\e \downarrow 0} \Epc \left[|u-u^\e|_{T, \infty}^p \right]=0, \ T>0,
\end{align*}
where $u(t,x)$ denotes the solution of the reflected SPDE \eqref{spde-1.1}, see Theorem 4.1 \cite{DoPa-93} or Theorem 2.1 \cite{XuZh-09}. In particular, we have
\begin{align}\label{eq-4.2-170611}
\lim_{\e \downarrow 0}\sup_{t\in [0,T]} \Epc \left[|u(t)-u^\e(t)|^2 \right]=0,\ T>0.
\end{align}
Therefore, to prove  Theorem \ref{thm-2.2-170519} and Theorem \ref{thm-2.3-170519}, it is important to show similar results hold for 
$P_t^\e\Phi$, where 
$$
P_t^\e\Phi(h)=\Epc \left[\Phi(u^\e(t; h)) \right], \ \Phi \in B_b(H).
$$ 
Let
$$F(u)(x)=f(u(x)), \  x\in [0,1], u\in H.$$ 
Then it is easy to know that 
\eqref{spde-approx} can be written in its abstract form.
\begin{align}\label{eq-3.2-1-170515}
\left\{
\begin{aligned}
du^\e(t) & =\left(Au^\e(t) + B(u^\e(t)) +\frac{1}{\e}F(u^\e(t) ) \right)dt +\Sigma(u^\e(t))dW(t),\\
u^\e(0)&=h. 
\end{aligned}
\right.
\end{align}

Let us now state the estimate of $|\nabla P_t^\e\Phi |^2$ for $t\geq 0$ and $\Phi \in C_b^1(H)$.
\begin{theorem} \label{thm-3.1-170514}
 Under the assumptions of {\bf A1)}-{\bf A3)},   we have that for all $\e>0$, the following holds.
\begin{align}\label{eq-3.5-1-170515}
|\nabla P_t^\e\Phi |^2 \leq 2M\exp\{\zeta (t) \} P_t^\e |\nabla\Phi|^2,\  \Phi\in C_b^1(H), t \geq 0,
\end{align}
where $M>0$ is the constant defined by \eqref{eq-3.16-170515} and $\zeta(t), t\geq 0$  is the function
defined  by \eqref{eq-2.3-170604}. 
\end{theorem}
\begin{proof}
The proof will be divided into two steps.\\
{\bf Step 1:} 
In this step, let us  further assume that $b$ and $\sigma$ are differentiable such that 
\begin{align}\label{eq-3.5-170514}
\sup_{u\in \R}|b'(u)|\leq L_b \ \ \mbox{and}\  \ \sup_{u\in \R}|\sigma'(u)|\leq L_\sigma.
\end{align}

Combining the assumptions {\bf A1)}-{\bf A3)} with the assumptions on $f$, we easily know that the operators $B, F$ and $\Sigma$ appeared in \eqref{eq-3.2-1-170515} satisfy all of the conditions in Theorem 5.4.1 {\rm (i)} in \cite{DaZa-96}. 
Hence, we have that $u^\e(\cdot; h)$ is  continuously differentiable in its initial data $h$.
Let $\nabla_k u^\e(t;h)$ denote the directional derivative of $u^\e(t;h)$ along the direction $k\in H$, that is, 
$$
\nabla_k u^\e(t;h) =\lim_{\delta \downarrow 0} 
\frac{u^\e(t; h +\delta k ) -u^\e(t;h)}{\delta}.
$$
We claim that  the following estimate holds for all $\e>0$ and $h\in H$:
\begin{align}\label{eq-3.4-1-170515}
\Epc[|\nabla_k u^\e(t;h)|^2] \leq M|k|^2\exp\{\zeta(t)\},\ t \geq 0,
\end{align}
where the constant 
$M>0$ and the function $\zeta(t), t\geq 0$ are  same as those stated in this theorem.

By Theorem 5.4.1 {\rm (i)} in \cite{DaZa-96}, we also know that its directional derivative $\nabla_k u^\e(t;h)$ along $k \in H$ is the unique solution of the SPDE
 \begin{align}\label{eq-3.6-170514}
\left\{
\begin{aligned}
d\nabla_k u^\e(t;h)  =& \big(A\nabla_k u^\e(t;h) +D B(u^\e(t); h)\cdot \nabla_k u^\e(t;h) \big)dt \\
& + \frac{1}{\e}DF(u^\e(t); h) \cdot \nabla_k u^\e(t;h)dt  \\ & +D\Sigma(u^\e(t); h)\cdot \nabla_k u^\e(t;h)dW(t), \ t \geq 0,\\
\nabla_k u^\e(0;h)= & k,
\end{aligned}
\right.
\end{align} 
where $D B, DF$ and $D\Sigma$ denote the Fr\'echet derivatives of $B, F$ and $\Sigma$.  Clearly, \eqref{eq-3.6-170514} can be rewritten as 
 \begin{align}\label{eq-3.7-170514}
\left\{
\begin{aligned}
& \frac{\partial}{\partial t} \nabla_k u^\e(t,x;h) \\
 = & \frac12 \frac{\partial^2}{\partial x^2}\nabla_k u^\e(t,x;h) +b'(u^\e(t,x);h) \nabla_k u^\e(t,x;h) 
 \\
& + \frac{1}{\e}f'(u^\e(t,x;h)) \nabla_k u^\e(t,x;h) \\
& +\sigma'(u^\e(t,x); h)\nabla_k u^\e(t,x;h)\dot{W}(t,x),\    t> 0, x\in (0,1),\\
 & \nabla_k u^\e(t,0;h)= \nabla_k u^\e(t,1;h)=0, \  t\geq 0, \\
& \nabla_k u^\e(0,x;h)=  k(x),\ x\in [0,1].
\end{aligned}
\right.
\end{align} 
Let us first assume $k$ is non-negative, i.e., $k \in K_0$. By the linearity of \eqref{eq-3.7-170514},  we know  that if  $k(x)\equiv 0, x\in [0,1]$, then $ \nabla_k u^\e(t,x;h) \equiv 0$ is the unique solution of \eqref{eq-3.7-170514}. Hence,
by using the comparison theorem of SPDEs,  see Theorem 2.1 \cite{DoPa-93} for example, we have
the  solution $\nabla_k u^\e(t,x;h)$ of \eqref{eq-3.7-170514}  is almost surely non-negative for any $k(x)\in K_0$, that is,  
\begin{align}\label{eq-3.8-170514}
\nabla_k u^\e(t,x;h)\geq 0, \ t \geq 0, x\in [0,1] \ a.s.
\end{align}
On the other hand, by the non-increasing property of $f$, we have that $$f'(u) \leq 0, \ u \in \R.$$ 
 Hence, by the non-negativity of the solution \eqref{eq-3.7-170514}, see  \eqref{eq-3.8-170514}, and using the comparison theorem of SPDEs again, we have that
\begin{align}\label{eq-3.9-170514}
0 \leq \nabla_k u^\e(t,x;h)\leq v^\e(t,x) \ a.s.,
\end{align}
where $v^\e(t,x)$ denotes the unique solution of 
\begin{align*}
\left\{
\begin{aligned}
\frac{\partial  v^\e}{\partial t}(t,x)  = & \frac{1}{2} \frac{\partial^2 v^\e(t,x) }{\partial x^2} +b'(u^\e(t,x;h)) v^\e(t,x) 
 \\
& +\sigma'(u^\e(t,x; h))v^\e(t,x)\dot{W}(t,x),\   t> 0, x\in (0,1), \\
v^\e(t,0;h) = & v^\e(t, 1;h), \ t\geq 0, \\ 
v^\e(0;h)= &  k(x),\ x\in [0,1].
\end{aligned}
\right.
\end{align*} 
In the following, we use the mild solution of $v^\e(t,x)$, that is,
\begin{align}\label{eq-3.11-170514}
v^\e(t,x) = &\int_0^1 g(t,x,y)k(y)dy +\int_0^t\int_0^1
g(t-s,x,y)b'(u^\e(s,y;h)) v^\e(s,y)dsdy  \notag \\
& + \int_0^t\int_0^1 g(t-s,x,y) \sigma'(u^\e(s,y;h)) v^\e(s,y)W(dsdy),\ t\geq 0 \ a.s. \notag
 \\
=: & I_1(t) +I_2(t) +I_3(t). 
\end{align}
Firstly, the property of $g(t,x,y)$ gives  easily that  \
\begin{align*}
|I_1(t)| \leq e^{-\frac{\pi^2 t}{2} }|k|,\ t \geq 0.
\end{align*}
Secondly, using Cauchy-Schwarz's inequality and  \eqref{eq-3.5-170514}, we can easily obtain that
\begin{align*}
\mathbb{E}[|I_2(t)|^2] \leq L_b^2t \int_0^t e^{-\pi^2 (t-s)} \mathbb{E}[|v^\e(s)|^2]ds,\ t \geq 0.
\end{align*}
Finally, It\^{o}'s isometric property gives that 
\begin{align*}
\mathbb{E}[| I_3(t) |^2] &=\int_0^t\int_0^1\int_0^1g^2(t-s,x,y)
\Epc[|\sigma'(u^\e(s,y); h) v^\e(s,y)|^2]dsdxdy \\
& \leq L_\sigma^2 \int_0^t\int_0^1\int_\R q^2(t-s,x-y)
\Epc[|v^\e(s,y)|^2]dxdsdy  \notag
\\
& = L_\sigma^2 \int_0^t\int_0^1 \frac{1}{2\sqrt{\pi (t-s)}}
\Epc[|v^\e(s,y)|^2]dsdy   \notag\\
&= \frac{L_\sigma^2}{2\sqrt{\pi }} 
 \int_0^t \frac{1}{\sqrt{t-s}} \Epc [|v^\e(s)|^2]ds,\ t \geq 0, \notag
\end{align*}
where $q(t,x)$ denotes the density of the Gaussian distribution with mean $0$ and variance $t$, the assumption \eqref{eq-3.5-170514}  and $g(t,x,y)\leq q(t,x-y)$ have been used for the second line. Therefore, by the above estimates and \eqref{eq-3.11-170514}, we obtain that 
\begin{align} \label{eq-3.12-170514}
\Epc [|v^\e(t)|^2]  \leq  & 3|k|^2 +3L_b^2t \int_0^t \mathbb{E}[|v^\e(s)|^2]ds \\
& + \frac{3L_\sigma^2}{2\sqrt{\pi }} 
 \int_0^t \frac{1}{\sqrt{t-s}} \Epc [|v^\e(s)|^2]ds,\ t \geq 0.  \notag
\end{align}
Let $\psi^\e(t)=\Epc [|v^\e(t)|^2], t\geq 0$. Using the  relation \eqref{eq-3.12-170514}, we have 
\begin{align*}
 \int_0^t \frac{1}{\sqrt{t-s}} \psi^\e(s)ds 
 \leq  & 6|k|^2t^{1/2} 
 +3L_b^2  \int_0^t 
\int_0^s  \frac{s}{\sqrt{t-s}} \psi^\e(\theta) d\theta  \\
&+   \frac{3L_\sigma^2}{2\sqrt{\pi }} 
\int_0^t   \int_0^s \frac{1}{\sqrt{t-s}\sqrt{s-\theta}} \psi^\e(\theta) d\theta  ds \\
\leq & 6|k|^2t^{1/2}  + \left\{ 6L_b^2t^{3/2}+ \frac{3L_\sigma^2\sqrt{\pi}}{2}  \right\}\int_0^t \psi^\e(\theta) d\theta,
\end{align*}
where $\int_\theta^t \frac{1}{\sqrt{t-s}\sqrt{s-\theta}}ds =\pi$ has been used for the second inequality.

Inserting this estimate to the last part of \eqref{eq-3.12-170514}, we see that 
\begin{align*} 
\psi^\e(t) \leq & \left\{ 3+ \frac{9L_\sigma^2}{\sqrt{\pi} } t^{1/2} \right\} |k|^2 \\
&+ \left\{\frac{9L_\sigma^4}{4}+ 3L_b^2t + \frac{9L_b^2 L_\sigma^2}{\sqrt{\pi} }t^{3/2} \right\} \int_0^t \psi^\e(s)ds \notag \\
:=& \alpha(t)|k|^2 + \beta(t) \int_0^t \psi^\e(s)ds, \ t\geq 0. \notag
\end{align*}
Note that $\alpha(t)$ and $\beta(t)$ do not depend on $\e$.
Using the generalized Gronwall inequality, see Lemma \ref{lem-3.2-170514} below, we deduce  that for all $\e>0$, 
\begin{align*}
\psi^\e(t) \leq \alpha(t)|k|^2 + \beta(t)\int_0^t \alpha(s) |k|^2  \exp \left\{ \int_s^t \beta(\theta)d\theta \right\}ds, \ t\geq 0.
\end{align*}
Since $\alpha(t)$ is increasing in $t\in [0, \infty)$ and  $\beta(t) \geq \frac{9L_\sigma^4}{4}$ for all $t\geq 0$, we have 
\begin{align*}
\psi^\e(t) & \leq \alpha(t)|k|^2 + \alpha(t) \beta(t) |k|^2 \int_0^t  \exp \left\{ \int_s^t \beta(\theta)d\theta \right\}ds \\
& \leq \alpha(t)|k|^2 +\frac{4}{9L_\sigma^4} \alpha(t)  \beta(t) |k|^2 
\int_0^t  \beta(s)\exp \left\{ \int_s^t \beta(\theta)d\theta \right\}ds  \notag \\
& = \alpha(t)|k|^2 +\frac{4}{9L_\sigma^4} \alpha(t) \beta(t) |k|^2  
\left[ \exp \left\{ \int_0^t \beta(\theta)d\theta \right\} -1\right]. \notag
\end{align*}
Using  the relation 
$\beta(t) \geq \frac{9L_\sigma^4}{4}$ for all $t\geq 0$
again, the above estimate gives  that 
\begin{align}\label{eq-3.15-170515}
\psi^\e(t) & \leq \frac{4 \alpha(t) \beta(t)}{9L_\sigma^4} |k|^2  
\exp \left\{ \int_0^t \beta(\theta)d\theta \right\}, \ t\geq 0. 
\end{align}
By the definitions of $\alpha(t)$ and $\beta(t)$, it is easy to check that  
\begin{align*}
\frac{4  \alpha(t) \beta(t)}{9L_\sigma^4} \leq M\left\{1+t^{1/2} +\frac{t}{2!} + \frac{t^{3/2}}{3!}
+ \frac{t^{2}}{4!} \right\}  \leq M\exp\{t^\frac{1}{2} \},\ t\geq 0;
\end{align*} 
recalling that $M$ is  the constant defined by \eqref{eq-3.16-170515}.

Combining this estimate with  \eqref{eq-3.15-170515}, we obtain that
\begin{align*}
\psi^\e(t) & \leq M|k|^2  
\exp \left\{t^{1/2}+ \int_0^t \beta(s)ds\right\}, \ t\geq 0.
\end{align*}
Therefore,  by the definition of $\beta(t)$ and recalling the definition of  $\zeta(t)$, see \eqref{eq-2.3-170604}, we easily obtain the relation
\begin{align*}
\Epc [|v^\e(t)|^2]\leq {M}|k|^2\exp \{\zeta(t)\},\ t \geq 0.
\end{align*}

As a consequence of the above estimate and the relation \eqref{eq-3.9-170514},  we have
for all $h\in H$ and $ k\in K_0$, 
\begin{align}\label{eq-3.18-170515}
\Epc [|\nabla_ku^\e(t;h)|^2]\leq M|k|^2\exp\{\zeta(t)\},\ t \geq 0.
\end{align}
From now on, let us consider the general $k\in H$. In this case, we have
$$
k(x)=k^+(x) -k^-(x), \ x\in [0,1],
$$
where $k^+(x)=\max\{k(x),0 \}$ and  $k^-(x)=\max\{-k(x),0\}$.
On the other hand, we also have that
$$
\nabla_ku^\e(t;h)= \nabla_{k^+}u^\e(t;h)-\nabla_{k^-}u^\e(t;h).
$$
Therefore,  using \eqref{eq-3.18-170515},
we have that 
\begin{align*}
\Epc [|\nabla_ku^\e(t;h)|^2] &\leq
2\{  \Epc[|\nabla_{k^+} u^\e(t;h)|^2] +\Epc [|\nabla_{k^-}u^\e(t;h)|^2]  \} \\
& \leq 2M(|k^+|^2+ |k^-|^2)\exp\{\zeta(t) \} \notag \\
& = 2M|k|^2\exp\{\zeta(t) \},\ t \geq 0, \notag
\end{align*}
which is the desired result \eqref{eq-3.4-1-170515}.
\\
\noindent {\bf Step 2:} Let us formulate the proof of \eqref{eq-3.5-1-170515} under the assumptions of this theorem. To do that, let $\phi\in C_0(\R)$ denote a symmetric and positive mollifier and set 
$$\phi_n(x)=n\phi(nx),\ n\in \N.$$
Let us now construct approximating sequences  $b_n, f_n$ and $\sigma_n, n\in \N$ for $b, f$ and $\sigma$ as below.
\begin{align*}
b_n(u) =b*\phi_n(u), \ f_n(u)=f*\phi_n(u) \ \mbox{and}
 \ \sigma_n(u) =\sigma*\phi_n(u),\ u\in\R,
\end{align*}
where $*$ denotes the convolution of two functions. From the properties of the mollifier, it is easy to know that $b_n, f_n$
 and $\sigma_n$ are smooth,  and  as $n \to \infty$,
$b_n, f_n$ and $\sigma_n$ converge to $b, f$ and $\sigma$ respectively. In addition, by Rademacher's Theorem, see Theorem 6 
in Section 5.8.3 \cite{Eva}, and then the assumptions 
{\bf A1)}-{\bf A2)}, we have that 
\begin{align}\label{eq-3.20-170515}
\sup_{n\in\N, u\in \R } |b_n'(u)|\leq L_b \ \ \mbox{and} \ 
\sup_{n\in\N, u\in \R } |\sigma_n'(u)|\leq L_\sigma.
\end{align}
Moreover, by the assumptions of $f$, we know that 
$$f_n'(u)\leq 0,\ u\in \R.$$ 
Therefore, $b_n, f_n$ and $\sigma_n$ satisfy all of the assumptions used in {\bf Step 1}.  Let 
$u^{\e, n}(t,x;h)$ denote the solution of SPDE$(b_n, f_n, \sigma_n)$ \eqref{spde-approx} and $\nabla_ku^{\e, n}(t,;h)$ denote its  
directional derivative along $k \in H$.  
By \eqref{eq-3.20-170515} and the estimate  \eqref{eq-3.4-1-170515} established in {\bf Step 1}, we have that 
for all $\e>0, n\in \N$ and $h\in H$, the following holds.
\begin{align}\label{eq-3.21-170515}
\Epc [|\nabla_k u^{\e,n}(t;h)|^2] \leq 2M|k|^2\exp\{\zeta(t)\},\ t \geq 0.
\end{align}
Owing to \eqref{eq-3.20-170515} and the proof of \eqref{eq-3.4-1-170515},  it is easy to know that $M$ and $\zeta(t), t\geq 0$ are independent of $\e$ and $n$ and  are same as those in \eqref{eq-3.4-1-170515}.
Set 
$$
P_t^{\e, n}\Phi(h)=\Epc[\Phi(u^{\e,n }(t;h))],\ \Phi\in B_b(H).
$$ 
By  Cauchy-Schwarz's inequality and \eqref{eq-3.21-170515},
we have 
\begin{align*}
| \nabla_k P_t^{\e, n}\Phi(h) |^2 &=\Epc [\langle \nabla \Phi(u^{\e,n }(t;h)),  \nabla_k u^{\e,n}(t;h) \rangle ] \\
& \leq   P_t^{\e, n}|\nabla_k \Phi|^2(h) \Epc \left[ |\nabla_k u^{\e,n}(t;h)|^2 \right]  \\
&\leq 
2M|k|^2\exp\{\zeta(t)\} P_t^{\e, n}|\nabla_k \Phi|^2(h).
\end{align*} 
Noting now that 
$$
|\nabla P_t^{\e, n}\Phi(h) | =\sup\{ 
  \nabla_k P_t^{\e, n} (h): \  k \in H\ \mbox{with}\  |k|=1
\},
$$
we easily obtain that 
\begin{align}\label{eq-4.21-170605}
|\nabla P_t^{\e, n}\Phi |^2 \leq 2M\exp\{\zeta (t) \} P_t^{\e,n} |\nabla\Phi|^2,\  \Phi\in C_b^1(H). 
\end{align}
On the other hand, we can easily show that 
\begin{align*}
\lim_{n\to \infty} \sup_{t\in [0,T]} \Epc [|u^{\e, n}(t; h) -u^{\e}(t; h)|^2]=0,
\end{align*}
see, for example, Theorem 7.1 \cite{DaZa-96} or Theorem A.1 \cite{PeZa-95}.
Consequently, we can obtain the desired result \eqref{eq-3.5-1-170515} by letting $n\to \infty$ in \eqref{eq-4.21-170605}.
\end{proof}

Let us formulate the following generalized Gronwall inequality,  which is used in the proof of Theorem \ref{thm-3.1-170514} and state its brief proof for the reader's convenience.
\begin{lemma}\label{lem-3.2-170514}
Suppose functions $\psi(t), \alpha(t), \beta(t)$ and $\gamma(t)$ are non-negative and locally integrable on $[0,\infty)$. If the function $\psi(t)$ satisfies the relation 
\begin{align} \label{eq-3.13-1-170514}
\psi(t) \leq \alpha(t) +\beta(t)\int_0^t \gamma(s) \psi(s)ds,\ t\geq 0,
\end{align}
then we have 
\begin{align*} 
\psi(t) \leq \alpha(t) +\beta(t)\int_0^t \alpha(s) \gamma(s)\exp \left\{\int_s^t \beta(\theta) \gamma(\theta)d\theta \right\} ds, \ t\geq 0.
\end{align*}
\end{lemma}
\begin{proof}
The proof is very easy. In fact, let 
\begin{align*}
\phi(s) = \exp \left\{-\int_0^s \beta(\theta) \gamma(\theta)d\theta \right\} \int_0^s \gamma(\theta) \psi(\theta)d\theta,\ s\geq 0.
\end{align*}
Noting that
\begin{align*}
\phi'(s) =\gamma (s)  \exp \left\{-\int_0^s \beta(\theta) \gamma(\theta)d\theta \right\}
\left\{\psi(s)-\beta(s)\int_0^s \gamma(\theta) \psi(\theta)d\theta \right\},
\end{align*}
by the relation \eqref{eq-3.13-1-170514}, we easily know that 
\begin{align*}
\phi(t)\leq \int_0^t  \alpha(s)\gamma(s) \exp \left\{-\int_0^s \beta(\theta) \gamma(\theta)d\theta \right\} ds.
\end{align*}
Hence,  by the definition of $\phi$ above, we have that
\begin{align*}
\int_0^t \gamma(s) \psi(s)ds & =\phi(t)  \exp \left\{ \int_0^t \beta(\theta) \gamma(\theta)d\theta \right\} \\
& \leq \int_0^t  \alpha(s)\gamma(s) \exp \left\{\int_s^t \beta(\theta) \gamma(\theta)d\theta \right\} ds.
\end{align*}
Consequently, inserting the above estimate into the right hand side of 
\eqref{eq-3.13-1-170514}, we can easily obtain the desired  relation.
\end{proof}

\begin{theorem}\label{thm-3.3-170517}
Suppose the assumptions {\bf A1)}-{\bf A3)} are satisfied. Then, the log-Harnack inequality holds for $P_t^\e$. More precisely, for any strictly positive $\Phi\in B_b(H)$,
\begin{align}\label{eq-3.25-170514}
P_t^\e \log \Phi(h_1) \leq \log P_t^\e\Phi(h_2) + \frac{M|h_1-h_2|^2}{k_1^2\int_0^t \exp\{-\zeta(s) \}ds}, \ h_1, h_2\in H, t>0, 
\end{align}
where $M$ and $\zeta(t), t\geq 0$ are same as those in Theorem \ref{thm-3.1-170514}.
\end{theorem}

\begin{proof}

As we saw in the proof of Theorem \ref{thm-3.1-170514}, without loss of 
generality, we assume $b, f$ and $\sigma$ are twice differentiable with bounded derivatives.
Let us first claim that for each $\Phi\in C_b^2(H)$ and $\e>0$
\begin{align}\label{eq-3.26-170514}
&P_s^\e \log P_{t-s}^\e\Phi(h)  \\ =
&\log P_{t}^\e\Phi(h)  -
\frac{1}{2} \int_0^s  P_{\theta }^\e| \Sigma(h)^* D \log P_{t-\theta}^\e\Phi(h)|^2d\theta,\ s\in [0,t].  \notag
\end{align} 
To prove this assertion, let  $A_n$ be the Yosida
approximation of $A$, that is, $$A_n=nA(n-A)^{-1}$$ and let $\Pi_n$ be the  
orthogonal projection from $H$ to ${\rm span} \{e_1, e_2, \cdots, e_n \}, n\in \N$.
Consider the following approximating stochastic partial differential equation for the penalized SPDE \eqref{eq-3.2-1-170515}:
\begin{align}\label{eq-3.27-1-170517}
\left\{
\begin{aligned}
du_n^\e(t) & =\Big(A_nu_n^\e(t) + B(u_n^\e(t)) +\frac{1}{\e}F(u_n^\e(t)) \Big)dt +\Sigma(u_n^\e(t))\Pi_n dW(t),\\
u^\e(0)&=h. 
\end{aligned}
\right.
\end{align}
It is known that for each $\e>0$ and $n\in \N$, this equation has a unique strong solution (see Section 7.4 \cite{DaZa-92}) and moreover,
\begin{align}\label{eq-3.27-170515}
\lim_{n\to \infty} \sup_{t\in [0,T]}\Epc \left[|u_n^\e(t) -u^\e(t)|^2 \right]=0.
\end{align}
Moreover,  the solution $(u_n^\e(t))_{t\geq 0}$ of \eqref{eq-3.27-1-170517} can be represented by   the stochastic equation
\begin{align*}
u_n^\e(t) =& h +\int_0^t \left( A_nu_n^\e(s) + B(u_n^\e(s)) +\frac{1}{\e}F(u_n^\e(s))  \right)ds \\
 & + \int_0^t \Sigma(u_n^\e(s))\Pi_n dW(s)\ a.s. \notag
\end{align*} 
Let us introduce the Markov semigroup
$$
P_{n, t}^\e\Phi(h) =\Epc \left[\Phi(u_n^\e(t; h) ) \right], \ \Phi\in B_b(H), t\geq 0.
$$
Then, by Theorem 5.4.2 \cite{DaZa-96}, for each $\Phi\in C_b^2(H)$, it is known that 
 $$P_{n, t}^\e\Phi(h) \in C_b^{1,2}([0,\infty) \times H).$$ Define the operator 
$\mathcal{L}$ on $C_b^2(H)$ by 
\begin{align*}
\mathcal{L}\Phi(h)  = &
\left\langle 
A_n(h)+ B(h) +\frac{1}{\e}F(h), D\Phi(h) 
\right\rangle 
\\
&+
\frac12{\rm Tr}\left[D^2\Phi(h) (\Sigma(h)\Pi_n )(\Sigma(h)\Pi_n )^* \right],\  h\in H.
\end{align*}
Then, it is easy to see that for any strictly positive and bounded function $\Phi$ on $H$,  
\begin{align*}
 \mathcal{L}\log P_{n,t}^\e \Phi (h)
=
\frac{\mathcal{L} \log P_{n,t}^\e  \Phi(h)}{\log P_{n,t}^\e  \Phi(h)}- \frac{1}{2}| (\Sigma(h)\Pi_n)^*  D\log P_{n,t}^\e  \Phi(h)|^2.
\end{align*} 
Hence, applying It\^{o}'s formula (for example, 
see Theorem 4.32 \cite{DaZa-92}) to the solution $u_n^\e(s)$ and 
$\log P_{n,t-s}^\e\Phi(h), s\leq t$, we have that
\begin{align}\label{eq-3.29-170515}
d\log P_{n,t-s}^\e \Phi(u_n^\e(s))  = &\partial_{s}\log P_{n,t-s}^\e  \Phi (u_n^\e(s))ds
+ \mathcal{L}\log P_{n,t-s}^\e  \Phi (u_n^\e(s))ds \\
& + \left\langle 
D\log P_{n,t-s}^\e  \Phi (u_n^\e(s)), \Sigma(u_n^\e(s))\Pi_n dW(s)
\right\rangle \notag\\
= & - \frac{1}{2}|( \Sigma(u_n^\e(s))\Pi_n)^* D\log P_{n,t-s}^\e  \Phi (u_n^\e(s))|^2 ds \notag \\
& +\left\langle 
D\log P_{n,t-s}^\e  \Phi (u_n^\e(s)), \Sigma(u_n^\e(s))\Pi_n dW(s)
\right\rangle. \notag
\end{align}
Noting the last term in \eqref{eq-3.29-170515} is a martingale, integrating  of both sides of \eqref{eq-3.29-170515} from $0$ to $s$ and then taking expectation,  we have that
\begin{align*}
& \Epc \left[\log P_{n,t-s}^\e(u_n^\e(s)) \right] 
\\
=& \log P_{n,t}^\e\Phi(h) 
-
\frac{1}{2} \int_0^s\Epc \left[ | (\Sigma(u_n^\e(\theta))\Pi_n)^* D\log P_{n,t-\theta }^\e \Phi(u_n^\e(\theta))|^2 \right]d\theta, \ s\in [0,t].
\end{align*} 
Using the notation of the Markov semigroup $P_{n,t}^\e$, we can rewrite the above equality as the following.
\begin{align*}
& P_{n,s}^\e \log P_{n,t-s}^\e(h) \\=
& \log P_{n,t}^\e\Phi(h) -
\frac{1}{2} \int_0^t P_{n,\theta }^\e| (\Sigma(h)\Pi_n)^* D\log P_{n,t-\theta}^\e\Phi(h)|^2d\theta,\ s\in [0,t]. \notag  
\end{align*} 
Recalling \eqref{eq-3.27-170515},  we can conclude the proof of \eqref{eq-3.26-170514} by letting $n\to \infty$ in the above equation.

From now on, let us formulate the proof of \eqref{eq-3.25-170514} based on the main ideas used in  Theorem  1.1 \cite{RoWa-10}.
Define $\alpha(s)$ by 
$$
\alpha(s) = \frac{\int_0^s \exp\{-\zeta(\theta)\} d\theta }{\int_0^t \exp\{-\zeta(\theta)\} d\theta},\ s\in [0,t].
$$
Then, it is easy to know that $\alpha(s)\in C^1([0,t])$ is increasing in $s$ and  $\alpha(0)=0,\ \alpha(t)=1$. For any $h_1, h_2\in H$, define
$h(s)\in H$ by 
$$
h(s)=h_1 \alpha(s) +(1-\alpha(s))h_2, \ s\in [0,t].
$$
It is clear that $h(0)=h_2$ and $h(t)=h_1$.
Using  \eqref{eq-3.26-170514} and the relation
\begin{align*}
P_t^\e \log\Phi(h_1) =\log P_t^\e\Phi(h_2) +
\int_0^t \frac{d}{ds} \left\{ \left( P_s^\e \log P_{t-s}^\e \Phi\right) 
(h(s) ) \right\}ds,
\end{align*}
we have that
\begin{align*}
& P_t^\e \log\Phi(h_1) - \log P_t^\e\Phi(h_2)\\
= & -\frac12\int_0^t P_{s}^\e| \Sigma(\cdot)^* D\log P_{t-s}^\e \Phi |^2(h(s))ds\\
&  + \int_0^t \alpha'(s) \langle DP_s^\e \log P_{t-s}^\e \Phi (h(s) ), h_1-h_2 \rangle ds. 
\end{align*}
Noting that the assumption $0<\kappa_1 \leq |\sigma(u)|$, see the assumption {\bf A3)}, we have $\kappa_1 \leq |\Sigma(h)|$ for all $ h\in H$. Therefore, by Theorem  \ref{thm-3.1-170514}, we obtain  that
\begin{align*}
& P_t^\e \log\Phi(h_1) - \log P_t^\e\Phi(h_2)\\
\leq  &
-\frac{\kappa_1^2}{2}\int_0^t P_{s}^\e|  D\log P_{t-s}^\e  \Phi |^2(h(s))ds\\
&  + \int_0^t \alpha'(s) |h_1-h_2 |
|DP_s^\e \log P_{t-s}^\e\Phi| (h(s) ) ds \\
\leq  &
-\frac{\kappa_1^2}{4}\int_0^t \frac{\exp\{-\zeta(s)\}}{M} |\nabla P_{s}^\e \log P_{t-s}^\e|^2(h(s))ds\\
&  + \int_0^t \alpha'(s) |h_1-h_2 |
|\nabla P_s^\e \log P_{t-s}^\e \Phi| (h(s) )ds \\
\leq & \frac{M|h_1-h_2|^2}{\kappa_1^{2} } \int_0^t\exp\{\zeta(s)\} (\alpha'(s))^2 ds\\
= &  \frac{M|h_1-h_2|^2}{\kappa_1^2\int_0^t \exp\{-\zeta(s) \}ds}, 
\end{align*}
where the Young inequality $|ab| \leq 2^{-1}\delta a^2 + {(2\delta)}^{-1} b^2,\  a, b\in \R, \ \delta >0$  and the definition of $\alpha(s)$ have been used for the third inequality and the last equality respectively.
Consequently, the proof of this theorem is completed.
\end{proof}

Based on Theorem \ref{thm-3.1-170514} and Theorem \ref{thm-3.3-170517}, we can easily describe the proofs of  Theorem 
\ref{thm-2.2-170519} and Theorem \ref{thm-2.3-170519}. 
\begin{proof} {(\bf Proofs of  Theorem 
\ref{thm-2.2-170519} and Theorem \ref{thm-2.3-170519}):}
Noting the relation \eqref{eq-4.2-170611} and letting $\e \downarrow 0$ in \eqref{eq-3.5-1-170515}, we can obtain \eqref{eq-2.4-170611} and then complete the proof of Theorem \ref{thm-2.2-170519}.

Let us turn to the proof of Theorem \ref{thm-2.3-170519}.
 By the relation \eqref{eq-4.2-170611} and letting $\e \downarrow 0$
 in \eqref{eq-3.25-170514}, we  can easily see that  \eqref{eq-2.4-170520} holds for any Lipschitz continuous function 
$\Phi$ and it can been extended for all $\Phi \in B_b(K_0)$ by the 
monotone class theorem. Thus, the proof of  Theorem 
\ref{thm-2.3-170519} is completed. 
\end{proof} 

\begin{remark} \label{rem-4.1-170607}
Noting that $f$ is Lipschitz continuous with the Lipschitz constant $1$, under the assumptions  {\bf A1)}-{\bf A3)}, we can easily testify  that for any $h_1, h_2 \in H$ and $t>0$, the following are fulfilled.
\begin{align*}
&\left|T_t\left(B(h_1)+\frac{1}{\e}F(h_1) \right)-T_t\left(B(h_2)+\frac{1}{\e}F(h_2)\right) \right|^2 \leq \left(L_b+\frac{1}{\e}\right)^2 |h_1- h_2|^2,\\
& \|T_t\Sigma(h_1) -T_t\Sigma(h_1)\|_{HS}^2 \leq 2L_\sigma^2 |h_1- h_2|^2\sum_{n=1}^\infty  \exp\{-n^2\pi^2t\}, 
\end{align*}
where $T_th(x) =\int_0^1 g(t,x,y)h(y) dy, h\in H$.
Therefore, for each $\e>0$, according to Theorem 1.1 and  Theorem 4.1 \cite{WaZh-14}, we can easily obtain that 
\begin{align}\label{eq-4.22-170607}
|\nabla P_t^\e\Phi |^2 \leq 6^{1+\frac{t}{t_0(\e)}} P_t^\e |\nabla\Phi|^2,\  \Phi\in C_b^1(H)
\end{align}
and the log-Harnack inequality holds for $P_t^\e$, i.e., for all $0\leq \Phi\in B_b(H)$,
\begin{align}\label{eq-4.23-170607}
& P_t^\e \log \Phi(h_1)  \\ \leq&  \log P_t^\e\Phi(h_2)
+ \frac{3 |h_1-h_2|^2 \log 6}{\kappa_1^2t_0(\e)\left( 1-6^{-t/t_0(\e)} \right) }, \ h_1, h_2\in H, t>0,  \notag
\end{align}
where \begin{align*}
& t_0(\e)\\
:=&  \sup\left\{t>0: \frac{\left(L_b+\frac{1}{\e} \right)^2t}{\pi^2} (1-\exp\{-\pi^2t\})
 +\frac{2L_\sigma^2}{\pi^2} \sum_{n=1}^\infty 
\frac{1-\exp\{-n^2\pi^2t\} }{n^2} \leq \frac16
\right\}.     
\end{align*}
However, it is clear that $t_0(\e)>0$ converges to $0$ as $\e \to 0$. 
Consequently, we can not deduce our many results in Theorem 
\ref{thm-2.2-170519} and Theorem \ref{thm-2.3-170519} 
 from  \eqref{eq-4.22-170607} and  \eqref{eq-4.23-170607}.
\end{remark}

\section{Proof of Theorem   \ref{thm-2.4-170605} }\label{sec-5}
As in Section \ref{sec-4}, let us first show the following results for the Markov semigroup $P_t^\e$ relative to the solution of the penalized SPDE \eqref{spde-approx}.
\begin{theorem} \label{thm-5.1-170605}
Suppose the assumptions {\bf A1)}-{\bf A3)} are satisfied. Then, we have
\\
{\rm (i)} For any $\Phi\in C_b^1(H)$,
\begin{align}\label{eq-3.31-170514}
P_t^\e\Phi^2 -(P_t^\e\Phi)^2\leq  2Mk_2^2 P_t^\e|\nabla\Phi|^2  \int_0^t \exp\{\zeta(s) \}ds, \ t>0.
\end{align}
{\rm(ii)} For any $\Phi\in B_b(H)$,
\begin{align}\label{eq-3.30-170514}
|\nabla P_t^\e\Phi|^2 \leq  \frac{ 2M \left(P_t^\e\Phi^2 -(P_t^\e\Phi)^2 \right)}{ k_1^2 \int_0^t \exp\{-\zeta(s) \}ds}, \ t>0.
\end{align}
\end{theorem}
\begin{proof}
The main idea to prove this theorem is similar to that used in 
Theorem \ref{thm-3.3-170517}.  We assume that the coefficients
 $b, f$ 
and $\sigma$ are twice continuously differentiable with bounded 
derivatives and use the same notations introduced in the proof of 
Theorem \ref{thm-3.3-170517}. Considering  the same approximating 
equation \eqref{eq-3.27-1-170517} as we did and then  applying It\^{o}'s formula, we have that
\begin{align*}
d(P_{n, t-s}^\e\Phi)^2(u_n^\e(s)) = & \partial_s (P_{n, t-s}^\e\Phi)^2(u_n^\e(s)) +\mathcal{L}(P_{n, t-s}^\e\Phi)^2(u_n^\e(s))\\
& + \langle D(P_{n, t-s}^\e\Phi)^2(u_n^\e(s)), \Sigma(u_n^\e(s))\Pi_n dW(s) \rangle\\
 = & |(\Sigma(u_n^\e(s))\Pi_n)^* DP_{n, t-s}^\e\Phi)(u_n^\e(s))|^2\\
& + \langle D(P_{n, t-s}^\e\Phi)^2(u_n^\e(s)), \Sigma(u_n^\e(s))\Pi_n dW(s) \rangle, s\in [0,t].
\end{align*}
Then, integrating both sides of the above equation from $0$ to $s\leq t$ and then taking expectation, we obtain that 
\begin{align*}
P_{n,s}^\e(P_{n, t-s}^\e\Phi)^2(h)= (P_{n,t}^\e\Phi)^2(h) + \int_0^s P_{n, \theta}|(\Sigma(\cdot)\Pi_n)^* DP_{n, t-\theta}^\e\Phi)(\cdot)|^2(h)d\theta.
\end{align*} 
In particular, let $s=t$ in the above equation and then let $n\to \infty$, we have
\begin{align}\label{eq-3.33-170517}
P_{t}^\e \Phi^2(h)= (P_{t}^\e \Phi)^2(h) + \int_0^t P_{ s}^\e|\Sigma(\cdot)^* \nabla P_{t-s}^\e\Phi)(\cdot)|^2(h)ds, \ t>0.
\end{align} 
Then we can easily complete our proofs. 
Let us first prove {\rm (i)} by using this relation.
In fact, since for all $u\in \R$, $|\sigma(u)|\leq \kappa_2$, from 
\eqref{eq-3.33-170517} and  Theorem  \ref{thm-3.1-170514}, it follows that 
\begin{align*}
P_{t}^\e  \Phi^2(h)& \leq  (P_{t}^\e \Phi)^2(h) + \kappa_2^2 \int_0^s P_{ s}^\e| \nabla P_{t-s}^\e\Phi|^2(h)ds\\
& \leq  (P_{t}^\e \Phi)^2(h) + 2M \kappa_2^2 \int_0^t P_{ s}^\e P_{t-s}^\e| \nabla \Phi|^2(h) \exp\{\zeta (t-s) \}ds \\
& =  (P_{t}^\e \Phi)^2(h) + 2M \kappa_2^2 P_{t}^\e| \nabla \Phi|^2(h) \int_0^t \exp\{\zeta(s) \}ds, \ t>0,
\end{align*}
which completes the proof of {\rm (i)}.

Finally, let us briefly formulate the proof of  {\rm (ii)}. 
By  \eqref{eq-3.5-1-170515} and the assumption {\bf A3)}, we have that 
\begin{align*}
P_{ s}^\e| \Sigma(\cdot)^* \nabla P_{t-s}^\e\Phi|^2(h)
 & \geq  \kappa_1^2 P_{s}^\e | DP_{t-s}^\e\Phi |^2(h) \\
& \geq \frac{\kappa_1^2 \exp\{-\zeta(s) \} }{2M} |\nabla P_t^\e \Phi|^2(h), \ s\in[0,t].
\end{align*}
Inserting this into the right hand side of  \eqref{eq-3.33-170517}, 
we have 
\begin{align*}
P_{t}^\e  \Phi^2(h)\geq  (P_{t}^\e \Phi)^2(h)
 + \frac{\kappa_1^2 \int_0^t \exp\{-\zeta(s)\} ds}{2M} |\nabla P_t^\e \Phi|^2(h),
\end{align*}
which gives the desired result \eqref{eq-3.30-170514}.
\end{proof}

In the end, let us conclude this paper with the proof of Theorem \ref{thm-2.4-170605}.
\begin{proof}
({\bf Proof of Theorem \ref{thm-2.4-170605}):} It is easy to show
this theorem by using Theorem \ref{thm-5.1-170605}. In fact, as the 
proof of Theorem \ref{thm-2.3-170519}, by  the monotone class 
theorem, it is enough to show  \eqref{eq-2.7-170611} holds for any 
Lipschitz continuous function $\Phi$ on $K_0$. However, noting 
\eqref{eq-4.2-170611}, they can be easily obtained by letting $\e 
\downarrow 0$ in  \eqref{eq-3.31-170514}. On the other hand, 
\eqref{eq-2.6-170611} can be obtained by letting  $\e \downarrow 
0$ in \eqref{eq-3.30-170514}. Thus the proof of this theorem is 
completed.
\end{proof}

\vskip 0.5cm
\begin{center}
\noindent {\bf Acknowledgements:} 
\end{center}
{\it The author was supported  in part by
Grant-in-Aid for Scientific Research (C) 16K05197 and 	20K03627 from Japan Society for the Promotion of Science(JSPS).
}

\bibliographystyle{plain}

\end{document}